\documentclass{article}
\usepackage{CJKutf8}

\usepackage{amscd,amssymb,bbm,graphicx,epsfig,psfrag,epic,eepic,latexsym,amsmath}
\usepackage{amsmath}
\usepackage[arrow,matrix,curve]{xy}
\usepackage{mathrsfs}
\usepackage[normalem]{ulem}
\usepackage{soul}
\usepackage
{color}

\newcommand{\R}{\mathbb{R}}
\newcommand{\C}{\mathbb{C}}

\begin{document}

\title{Integrability of Kepler Billiards at Zero-Energy}
\author{ Lei Zhao}

\maketitle

\begin{abstract}
	We consider a Kepler billiard with zero-energy in the plane defined inside a smooth closed connected simple curve which intersects all focused parabola at at most two points.  {We show that} if has an invariant curve consisting of $2$-periodic orbits and there exists a $C^{1}$-first integral with non-vanishing gradient in the region between the invariant curve and the boundary curve, then the system is defined actually inside an ellipse with the Kepler center occupying one of the foci. This statement is obtained as a simple ``translation'' of the theorem of Bialy-Mironov \cite{BM} with Levi-Civita transformation. 
\end{abstract}


A Kepler billiard in the plane is the combination of a Kepler problem together with a piecewise-smooth curve, called the wall in the plane. The moving particle moves under the attraction of the Kepler center and gets reflected elastically when meeting the wall. 

The Hooke-Kepler correspondence, obtained e.g. via the Levi-Civita regularization \cite{LC} mapping has been noticed to extends to Kepler billiards in \cite{P}. Indeed the complex square mapping sends a pair of centrally-symmetric orbits of the Hooke problem to an orbit of the Kepler problem. {Moreover, as the complex square mapping is conformal, it preserves angles. Thus the law of elastic reflection defined by the equality of the angles of incidence of reflection is preserved as well.} Thus a Hooke billiard system in a centrally-symmetric domain is sent to a Kepler billiard system, with a pair of centrally-symmetric billiard orbits sent to the same billiard orbit in the Kepler billiard. Vice versa, the pre-image of each Kepler billiard orbit is a pair of centrally-symmetric Hooke billiard orbits. Here the central-symmetry is defined with respect to the Hooke center. 

In \cite{TZ} we conjectured that if a Kepler billiard in the plane bounded by a  smooth connected closed simple curve is integrable, then the bounding curve is an ellipse and the Kepler center is at one of the foci of this ellipse.  In \cite{BCB}, this integrable rigidity result has been established for Kepler billiard inside an ellipse with sufficiently large energy. Recently this has been extended to closed centrally-symmetric or cyclic-rotational-symmetric convex domains bounded by analytic curves with sufficiently large energy \cite{BBBT}. Here we show that the conjecture holds at zero energy among a class of domains that we call $K$-convex domains. 

To start with, we consider a Kepler billiard in an elliptic domain in $\R^{2}$. We identify $\R^{2} \cong \C$ and assume that the Kepler center is at $O$. By normalization, we assume that the underlying Kepler problem is given by the Hamiltonian  
$$H(p, q):=\dfrac{|p|^{2}}{2} - \dfrac{1}{|q|},$$
{defined on the phase space $\C \setminus O \times \C$ equipped with the standard symplectic form $\Re(d \bar{p} \wedge d q)$.}

We fix its energy to zero and we apply Levi-Civita transformation on this energy level. This means to change time on the zero energy hypersurface  $\{H(p, q)=0\}$ {(on which the orbits are focused parabolae)} by multiplying the Hamiltonian by the factor $|q|$, then pull the resulting Hamiltonian $|q| H (p, q)$ back by the canonical mapping 
$$\C \setminus O \times \C \to \C \setminus O \times \C, (z, w) \mapsto (q=z^{2}, p=w/(2 \bar{z})).$$

{Indeed it is now direct to check that $\Re(\bar{p} d q)=\Re(\bar{w} d z)$ holds and consequently $\Re(d \bar{p} \wedge d q)=\Re(d \bar{w} \wedge d z)$.}

{The transformed Hamiltonian just $\dfrac{|w|^{2}}{8}$ at its $1$-energy level, whose orbits are straight lines. This way we connect straight lines to focused parabolae via the complex square mapping.}

{Next we see how this transformation connects billiard systems:}


{The pre-image of any billiard table under the mapping $\C \to \C, z \mapsto z^{2}$ is a centrally-symmetric domain. As the complex square mapping sends a centered ellipse into a focused ellipse, If the Kepler billiard is inside an elliptic domain that is not focused at O, then the pre-image is not an elliptic domain.}

{A billiard trajectory of the Kepler billiard is formed by a sequence of Kepler parabolic arcs $\mathcal{K}_{i}$, concatenated by the condition that the end point of $\mathcal{K}_{i}$ agrees with the starting point of $\mathcal{K}_{i+1}$, this point lies in the wall of reflection and the angle of incidence agrees with the angle of reflection. The pre-image under the complex square mapping of each $\mathcal{K}_{i}$ are two centrally symmetric line segments $\mathcal{H}^{\pm}_{i}$ having start and end points lie centrally-symmetrically in the corresponding wall of reflection, which is the pre-image of the Kepler wall of reflection. As the domain contains $O$, the line segments $\mathcal{H}^{\pm}_{i}$ are disjoint. Choose one of these two line segments and call it $\mathcal{H}_{i}$. Then for $\mathcal{H}_{i+1}$ we choose the one out of $\mathcal{H}^{\pm}_{i+1}$ whose start point is the end point of $\mathcal{H}_{i}$. As the complex square mapping is conformal, the angle of incidence agrees with the angle of reflection at this point and thus $\mathcal{H}_{i}. \mathcal{H}_{i+1}$ is part of a billiard trajectory. This way we get two billiard trajectories from the pre-image of a Kepler parabolic billiard trajectory. }

In this way we have tranformed the zero-energy Kepler billiard into the usual Birkhoff billiard in a centrally-symmetric domain in the plane, bounded by a centrally-symmetric $C^{\infty}$-curve, up to a $2:1$-cover. 

Moreover, if the Kepler billiard admits a $C^{1}$-first integral {in a region near the boundary}, then the pull-back of this function is a $C^{1}$-first integral of the resulting Birkhoff billiard {in the corresponding region near the boundary.}

A {celebrated} theorem of Bialy-Mironov \cite{BM} asserts that if the Birkhoff billiard  in a centrally-symmetric domain in the plane has an invariant curve of $4$-periodic points, and the region $\mathcal{A}$ between this invariant curve and the boundary admits a $C^{1}$-first integral with non-vanishing gradient on $\mathcal{A}$, then the domain is an ellipse.
{It has been shown in \cite[Thm 4.1]{BM} that $4$-periodic orbits in a centrally-symmetric domain are centrally-symmetric parallelograms, which correspond to $2$-periodic orbits in the corresponding zero-energy Kepler billiard. So the invariant (caustic) curve of $4$-periodic orbits corresponds to the invariant curve of $2$-periodic orbits of the corresponding zero-energy Kepler billiard. A $C^{1}$-first integral with non-vanishing gradient on the region bounded by the invariant curve of $2$-periodic orbits to the boundary is pulled back to a (centrally-symmetric) $C^{1}$-first integral with non-vanishing gradient on $\mathcal{A}$. }

\medskip

From these considerations we deduce

\medskip
{\bf Theorem A:} {If a Kepler billiard in an elliptic table in $\R^{2}$ has an invariant curve of $2$-periodic points and the region bounded by this invariant curve to the boundary admits an additional $C^{1}$-first integral with non-vanishing gradient on its zero-enegy level, then the elliptic table is focused at the Kepler center.}
\medskip

The same is true for a more general type of domains. We call a closed domain in $\R^{2}$ having the Kepler center $O$ in its interior \emph{K-convex}, if the zero-energy Kepler parabolic orbits intersects its boundary at at most two points. This is a natural convexity assumption with respect to parabolic Keplerian orbits. For example, an elliptic domain containing $O$ at its focus is K-convex. {As a parabola may intersect an ellipse at 4 points, there are strictly-convex domains that are not $K$-convex. For example, fixing a parabola focused at $O$ and consider a very eccentric elliptic domain whose minor axis agrees with the axis of symmetry of the parabola. By increasing the semi major axis length of the elliptic domain and decrease the semi minor axis length we can arrive eventually at an elliptic table which is strictly-convex but not $K$-convex. It seems tricky to find bounded domains which are $K$-convex but are not strictly-convex. Should we release the boundedness assumption, then a half-plane containing $O$ is such an example.  }

{The pre-image of a K-convex domain under the complex square mapping is strictly convex and centrally-symmetric. Vice versa, the image of a strictly convex and centrally-symmetric domain containing $O$ under the complex square mapping is $K$-convex. As convexity is a $C^{2}$-open condition, being $K$-convex is also $C^{2}$-open. So a $C^{2}$-small deformation of an elliptic domain with $O$ at a focus is again \emph{K-convex} and thus there are many $K$-convex domains.}

Applying the same argument as above, we have

\medskip
{\bf Theorem B:} {If a Kepler billiard in a K-convex table in $\R^{2}$ has an invariant curve of $2$-periodic orbits and the region between this invariant curve and the boundary admits an additional $C^{1}$-first integral on its zero-enegy level with non-vanishing gradient, then the table is bounded by an ellipse focused at the Kepler center.}
\medskip

The notion of $K$-convexity domains extend directly to positive energies for Kepler problem in the plane, defined now with respect to hyperbolic orbits. For negative energies it is less clear as to obtain natural extensions of $K$-convexity we may need to allow unbounded domains, for which the analysis will be of very different nature. Nevertheless we expect that variations of Theorem B hold in these cases as well. 

\medskip

{\bf Acknowledgement:} L.Z. is supported by DFG ZH 605/4-1 and Fundamental Research Funds for the Central Universities of China.



School of Mathematical Sciences, Dalian University of Technology, China, 

Email: lei.zhao@math.uni-augsburg.de; zhao1899@dlut.edu.cn

\end{document}